\documentclass[12pt,a4paper]{article}
\usepackage{theorem}
\newtheorem{lemma}{Lemma}
\newtheorem{proposition}[lemma]{Proposition}
\newtheorem{theorem}[lemma]{Theorem}
\newtheorem{corollary}[lemma]{Corollary}
{\theorembodyfont{\upshape}}
{\theorembodyfont{\upshape}\newtheorem{remark}[lemma]{Remark}}
{\theorembodyfont{\upshape}}
{\theorembodyfont{\upshape}}

\newcommand{\Z}{{\bf Z}}
\newcommand{\R}{{\bf R}}
\newcommand{\C}{{\bf C}}
\newcommand{\rme}{{\rm e}}
\newcommand{\rmd}{{\rm d}}

\newcommand{\alp}{\alpha}
\newcommand{\bet}{\beta}

\newcommand{\lam}{\lambda}
\newcommand{\del}{\delta}
\newcommand{\eps}{\varepsilon}

\newcommand{\Dom}{{\rm Dom}}

\newcommand{\Spec}{{\rm Spec}}

\newcommand{\Nhd}{{\rm Nhd}}
\newcommand{\norm}{\Vert}
\renewcommand{\Re}{{\rm Re}\;}
\renewcommand{\Im}{{\rm Im}\;}

\newcommand{\Proof}{\underbar{Proof}{\hskip 0.1in}}

\newcommand{\dist}{{\rm dist}}
\newcommand{\Schrodinger}{Schr\"odinger }

\newcommand{\abs}[1]{\left\vert#1\right\vert}

\newcommand{\Arg}{{\rm Arg}\;}

\title{SPECTRAL PROPERTIES OF NON-SELF-ADJOINT OPERATORS IN
 THE SEMI-CLASSICAL REGIME}
\author{Paul Redparth}
\date{March 2000}
\begin{document}
\maketitle
\begin{abstract}
We give a spectral description of the semi-classical \Schrodinger
operator with a piecewise linear, complex valued potential.
Moreover, using these results, we show how an arbitrarily small
bounded perturbation of a non-self-adjoint operator can completely
change the spectrum of the operator.

\par
keywords: \Schrodinger operator, non-self-adjoint operator, Airy
function, characteristic determinant, asymptotic analysis.
\end{abstract}


\section{Introduction}\label{int}

This work was motivated by the paper of Shkalikov \cite{Shk97}
concerning the analysis of the semi-classical Airy operator,
together with our computer simulations of the associated discrete
problem using the numerical package Matlab. Specifically, we
examined the operator $H_{\del,h}$ given formally by
\begin{equation}\label{jumpop}
H_{\del,h}:=-h^2\frac{\rmd^2}{\rmd x^2}+V_{\del}\qquad\mbox{ on
}L^2(-1,1)
\end{equation}
where
\[ V_{\del}(x):=\left\{\begin{array}{ll} i(x+\del) &\mbox{ for $x>0$ }\\
i(x-\del) &\mbox{ for $x<0$ }\\ \end{array} \right. \] and both
$\del>0$, $h>0$ are small. In \cite{Shk97} it is shown rigorously
that as $h\to 0$ the spectrum of $H_{0,h}$ becomes dense inside an
arbitrarily small neighbourhood of the Y-shaped subset of $\C$
defined by
\[
[i,1/\sqrt{3}],[-i,1/\sqrt{3}]\qquad\mbox{ and  }\qquad
[1/\sqrt{3},\infty),
\]
where we use $[\alp,\bet]$ to denote the line segment joining
$\alp,\bet\in\C$ (see figure 2).

This paper confirms our surprising numerical results, which
suggested that when an arbitrarily small jump discontinuity is
inserted into the otherwise linear (purely imaginary) potential,
and then the semi-classical parameter $h$ allowed to go to $0$,
the asymptotic spectrum of $H_{\del,h}$ turns out to be completely
different from that of $H_{0,h}$ (see figure 1 and Corollary
\ref{c1}).

Several papers \cite{SHD91,Reddy93,Shk97} have been written about
the operator $H_{0,h}$ - a major motivation being that it is a
model operator for the `Orr-Sommerfeld' problem \cite{Reddy93}.
The operator also defines the `Squire model for the Couette flow'
in hydrodynamics; and in its own right, defines the semigroup
which is the solution of the so-called `Torrey equation'
\cite{SHD91}, related to the diffusion of magnetic fields. Thus,
although the spectrum of this non-self-adjoint operator displays
sometimes strange and singular behaviour, it must not be dismissed
as a `pathological' example from pure mathematics since it has
important applications - for example, in magnetic resonance
imaging devices.

It is well known that a basis for solutions of the
so-called `Airy equation'
\[ -f''(z)+zf(z)=0 \]
is given by any two of the Airy functions $Ai(z)$, $Ai(\rme^{-2\pi
i/3}z)$ and $Ai(\rme^{2\pi i/3}z)$ (see \cite{Olver74}). Thus, the
analytic investigation of the eigenvalues will involve examining
the asymptotic behaviour of certain Airy functions, where the
expressions we use are WKB-type approximations. We will show that
the eigenvalues lie inside a certain subset of the complex plane,
which is intimately related to the $Stokes'$ $lines$ (or
$principal$ $curves$) of the problem (see \cite{Olver74}). The
proof will depend upon showing that for all $\lam$ outside this
subset, the eigenvalue problem
\[ H_hf(x)=\lam f(x) \]
has a well-defined Green's function. As in \cite{Shk97}, our
analysis uses the concept of the $characteristic$ $determinant$
which we will describe in Section \ref{s1}.

In Section \ref{s2} we analyse the spectral behaviour in the
semi-classical limit $h\to 0$ for a general complex-valued,
piecewise linear potential. The surprising result is that each
linear segment of the potential gives rise to a characteristic
`Y-shaped' set of eigenvalues; and the spectrum of the operator is
contained within the superposition of these `Y-shaped' sets. From
the point of view of applications, this means that whilst the
asymptotic spectrum is theoretically computable for an idealised
linear potential; in practice, any $arbitrarily$ small deviation
from the ideal can $completely$ change the spectrum. One way of
expressing this for an operator $H$ is in terms of the
$pseudospectral$ $sets$:
\[
\Spec_{\epsilon}(H):=\Spec(H)\cup\{z\in\C:\norm(H-z)^{-1}\norm
\geq\epsilon^{-1}\},
\]
i.e. the contour sets of the resolvent norm, with the convention
that $z\in\Spec(H)$ implies
\[ \norm(H-z)^{-1}\norm
:=\infty. \]
For many non-self-adjoint operators it has been demonstrated
\cite{Dav97a,Reddy93,Tre97} that the pseudospectral sets become
very large as some parameter varies, even though $z$ may be far
from the spectrum of the operator. This is equivalent to saying
that the spectrum is computationally unstable. Our aim in this
paper is to demonstrate for a relatively transparent case, the
mechanism behind this phenomenon; we believe the results to be
capable of extension to a more general class of piecewise analytic
potentials.

In Section \ref{s3} we provide an analysis of the simultaneous
limit as $h\to 0$ and $\del\to 0$ together.


\section{The Characteristic Determinant}\label{s1}

In this section we describe the characteristic determinant of the
operator $H_h$ defined by
\[  H_hf(x):=-h^2\frac{\rmd^2f(x)}{\rmd x^2}+V(x)f(x) \]
acting on $L^2(-1,1)$ with Dirichlet boundary conditions, $h>0$
small, and $V(x)$ the complex valued, $n$-times piecewise linear
function
\[ V(x):=\left\{\begin{array}{ll} m_1x+l_1 & x_0\leq x<x_1 \\
m_2x+l_2 & x_1<x<x_2 \\ \qquad \vdots & \qquad \vdots
\\ m_nx+l_n & x_{n-1}<x\leq x_n \\
\end{array} \right. \]
with $-1=x_0<x_1<\ldots <x_n=1$, and the $m_i$, $l_i$
$i=1,\ldots,n$ complex constants. The domain of the operator is
given precisely by
\begin{equation}\label{boundconds1}
\Dom(H_h)=\{f\in C[-1,1]:f(-1)=f(1)=0,f'\in C[-1,1] \mbox{ and
}f''\in L^2(-1,1)\}
\end{equation}
where the primes $'$ denote differentiation with respect to $x$,
and $f''$ is initially to be interpreted in the distributional
sense. A direct substitution shows that a basis of solutions for
the differential equation
\[
-h^2 f''(x)+(V(x)-\lam)f(x)=0
\]
where $V(x):=mx+l$; and $l$, $m$ are complex constants, is given
by any two of the Airy functions $Ai(w)$ and $Ai(\rme^{\pm 2\pi
i/3}w)$, where
\[
 w:=h^{-2/3}m^{-2/3}(V(x)-\lam).\\
\]
It follows that, in order to construct an eigenfunction of the
operator $H_h$, we seek constants $\alp_{i1},\alp_{i2}$
$i=1,\ldots,n$ such that the function
\[ f(x):=\left\{\begin{array}{ll} \alp_{11}u_{11}(x)+
\alp_{12}u_{12}(x)& x_0\leq x<x_1\\
\alp_{21}u_{21}(x)+\alp_{22}u_{22}(x)& x_1<x<x_2 \\ \qquad
\qquad\vdots & \qquad\vdots \\
\alp_{n1}u_{n1}(x)+\alp_{n2}u_{n2}(x) & x_{n-1}<x\leq x_n\\
\end{array} \right.
\]
satisfies all of the domain conditions (\ref{boundconds1}), where

\begin{equation}\label{defu}
u_{i1}(x):=Ai(\rme^{-2k\pi
i/3}h^{-2/3}m_i^{-2/3}((m_ix+l_i)-\lam))
\end{equation}
with $k\in\{-1,0,1\}$. For each $i=1,\ldots,n$, the functions
$u_{i2}$ are defined similarly, except that a different choice of
$k$ is to be taken from $\{-1,0,1\}$.

In addition to satisfying the boundary conditions $f(-1)=f(1)=0$,
$f$ must also be continuously differentiable, even at the points
$x_i$. From the power series definition \cite[p.54]{Olver74}, it
is clear that the Airy functions $Ai$ are analytic on the whole of
$\C$, and so the requirement that $f$ be continuously
differentiable reduces to the $2(n-1)$ simultaneous `matching'
conditions
\[
\alp_{i1}u_{i1}(x_i-)+\alp_{i2}u_{i2}(x_i-)-\alp_{(i+1)1}u_{(i+1)1}(x_i+)-
\alp_{(i+1)2}u_{(i+1)2}(x_i+)=0
\]
and
\[
\alp_{i1}u'_{i1}(x_i-)+\alp_{i2}u'_{i2}(x_i-)-\alp_{(i+1)1}u'_{(i+1)1}(x_i+)-
\alp_{(i+1)2}u'_{(i+1)2}(x_i+)=0.
\]
The boundary conditions $f(-1)=f(1)=0$ demand that
\[
\alp_{11}u_{11}(-1)+\alp_{12}u_{12}(-1)=0
\]
and
\[
\alp_{n1}u_{n1}(1)+\alp_{n2}u_{n2}(1)=0.
\]
Thus finding a solution of the differential equation
\[
-h^2 f''(x)+V(x)f(x)=\lam f(x)
\]
which satisfies all of the domain conditions (\ref{boundconds1}),
involves solving the matrix equation
\[
\pmatrix{ u_{11}(-1) & u_{12}(-1) & 0 & 0 & 0 & \cdots & 0 \cr
u_{11}(x_1) & u_{12}(x_1) & -u_{21}(x_1) & -u_{22}(x_1) & 0 &
\cdots & 0 \cr u'_{11}(x_1) & u'_{12}(x_1) & -u'_{21}(x_1)
&-u'_{22}(x_1) & 0 & \cdots & 0 \cr \vdots & \vdots & \ddots &
\ddots & \ddots & \vdots & \vdots \cr 0 & \cdots & 0 &
u_{(n-1)1}(x_{n-1}) & u_{(n-1)2}(x_{n-1}) & -u_{n1}(x_{n-1}) &
-u_{n2}(x_{n-1}) \cr 0 & \cdots & 0 & u'_{(n-1)1}(x_{n-1}) &
u'_{(n-1)2}(x_{n-1}) & -u'_{n1}(x_{n-1}) & -u'_{n2}(x_{n-1}) \cr  0 &
\cdots & 0 & 0 & 0 & u_{n1}(1) & u_{n2}(1) \cr } \]
\begin{equation}\label{matrix1}
\times\pmatrix
{\alp_{11} \cr \alp_{12} \cr \alp_{21} \cr \alp_{22} \cr
\vdots \cr \alp_{n1} \cr \alp_{n2} \cr }=
\pmatrix{ 0 \cr 0 \cr 0 \cr 0 \cr \vdots \cr 0 \cr 0 \cr }
\end{equation}
for the constants $\alp_{i1},\alp_{i2}.$ Note that we are taking
the left- and right-hand limits at the points $x_i$, although here
and subsequently we will abuse the notation in order to add
clarity, and simply write $u_{i1}(x_i)$ etc. It is the determinant
of the matrix on the left-hand side of (\ref{matrix1}) that we
shall call the characteristic determinant of the eigenvalue
problem defined by $H_h$.


\section{Airy functions and Stokes' lines}\label{s4}

For our proofs in the next section, we will need some notation and
basic properties of the Airy functions. In all that follows we let
the argument function $\Arg$ take principal values i.e.
\[ \Arg:\C\to(-\pi,\pi]. \]
If $\Arg(\bet-\alp):=\theta$, the subsets $Y(\alp,\bet)$ of $\C$
are to be constructed as follows: take the lines
\[ \alp+r\rme^{2\theta i/3}\qquad\mbox{ and }\qquad
\left\{\begin{array}{ll} \bet+r\rme^{2\theta i/3+2\pi i/3} &\mbox{
for }\theta<0 \\ \bet+r\rme^{2\theta i/3-2\pi i/3} &\mbox{ for
}\theta\geq 0 \\
\end{array} \right.
\]
as $r$ ranges in $[0,\infty)$, to their point of intersection,
$\Gamma$ say. Then, from $\Gamma$ extend the infinite line defined
by the set of $\lam\in\C$ such that
\begin{equation}\label{equals}
\Re\{(\rme^{-2\theta i/3}(\alp-\lam))^{3/2}\}
=\Re\{(\rme^{-2\theta i/3}(\bet-\lam))^{3/2}\}.
\end{equation}
The motivation for this set will become clear during our proofs;
in fact, it will be seen to comprise a curve asymptotic to the
line
\[
\left\{z\in\C:\Im(z)=\frac{\Im(\alp)+\Im(\bet)}{2}\right\}.
\]
Note for now, however, that (\ref{equals}) is $h$-independent. The
$\eps$-neighbourhood of any subset $T$ of $\C$ will be defined by
\[
\Nhd(T;\eps):=\{t+z:t\in T\mbox{ and }\abs{z}<\eps\}.
\]
We will use the well-known \cite{Olver74} asymptotic expansion of
the Airy function $Ai$, giving the WKB-type approximation:
\begin{equation}\label{asymptotics1}
Ai(z)=\frac{z^{-1/4}}{2\sqrt\pi}\exp\left(-\frac{2}{3}z^{3/2}\right)
(1+O(z^{-3/2}))
\end{equation}
as $\abs{z}\to\infty$, valid for $\abs{\Arg(z)}<\pi$; and where
the principal value of $z^{3/2}$ is taken. Following the notation
of Olver (see \cite[p.413]{Olver74}), we define
\[ S_0:=\{z:\abs{\Arg(z)}<\pi/3\} \]
\[ S_1:=\{z:\pi/3<\Arg(z)<\pi\} \]
\[ S_{-1}:=\{z:-\pi/3>\Arg(z)>-\pi\} \]
(suffixes enumerated modulo 3). One can check that for all complex
$z$ (and taking principal values), we have
\begin{equation}\label{id1}
\Re\{(\rme^{-2\pi i/3}z)^{3/2}\}=\left\{\begin{array}{ll}
-\Re\{(\rme^{2\pi i/3}z)^{3/2}\} &\mbox{ for }z\in S_1\cup
S_{-1}\\ \Re\{(\rme^{2\pi i/3}z)^{3/2}\} &\mbox{ for }z\in S_0 \\
\end{array} \right.
\end{equation}
\begin{equation}\label{id2}
\Re\{(\rme^{-2\pi i/3}z)^{3/2}\}=\left\{\begin{array}{ll}
-\Re\{(z)^{3/2}\} &\mbox{ for }z\in S_1\cup S_0 \\
\Re\{(z)^{3/2}\} &\mbox{ for }z\in S_{-1} \\ \end{array} \right.
\end{equation}
and
\begin{equation}\label{id3}
\Re\{(\rme^{2\pi i/3}z)^{3/2}\}=\left\{\begin{array}{ll}
-\Re\{(z)^{3/2}\} &\mbox{ for }z\in S_0\cup S_{-1} \\
\Re\{(z)^{3/2}\} &\mbox{ for }z\in S_1. \\ \end{array} \right.
\end{equation}
Then, putting
\begin{equation}\label{defaik}
Ai_{k}(z):=Ai(\rme^{-2k\pi i/3}z)
\end{equation}
the asymptotics (\ref{asymptotics1}) show that as
$\abs{z}\to\infty$, $\abs{Ai_k(z)}$ decreases exponentially for
$z\in S_k$, and increases exponentially for $z\in S_{k-1}\cup
S_{k+1}$. The boundaries of the sectors $S_k$ i.e. the rays
$t\rme^{\pm \pi i/3}$ and $t\rme^{\pi i}$ for $t\in[0,\infty)$,
are known as the Stokes' lines (or principal curves) of the
problem \cite[p.503]{Olver74}. Indeed, for the Airy equation
\[
-f''(z)+zf(z)=0
\]
the Stokes' lines are defined to be the values of $z$ such that
\[
\Re\int_0^z\sqrt{t}\;\rmd
t=\Re\left\{\frac{2}{3}z^{3/2}\right\}=0,
\]
and denote the boundaries of the $principal$ $subdomains$ $S_1$
etc., inside of which the asymptotic expression
(\ref{asymptotics1}) is valid for each $k$.

We will call the suffix $k$ `allowable' for any given $z\in\C$, if
\begin{equation}\label{allowable}
\abs{\Arg(\rme^{-2k\pi i/3}z)}<\pi.
\end{equation}


\section{The Semi-Classical Limit}\label{s2}

The following is a generalisation of the argument used in
\cite{Shk97} for the potential $V(x)=ix$, and will form a lemma
for our main theorem.

\begin{proposition}\label{t2}
Let $V$ be the complex valued linear potential given by
\[ V(x)=mx+l \qquad x\in[-1,1] \]
where $m$ and $l$ are complex constants; $u_{11}(x)$, $u_{12}(x)$
are as defined in (\ref{defu}), and $a,b\in [-1,1]$, $a<b$.
\par
Let $\eps>0$ be given and $\lam\in\C$. If
\[
\lam\notin\Nhd(Y(V(a),V(b));\eps),
\]
then
\begin{equation}\label{state1}
u_{11}(b)u_{12}(a)=o(u_{11}(a)u_{12}(b))
\end{equation}
as $h\to 0$.
\end{proposition}
\Proof A simple scaling and translation of the operator $H_h$
allows us, without loss of generality, to assume that $a:=-1$,
$b:=1$ and $l=0$. That is, we assume
\[
V(x):=x\rme^{i\theta},\qquad\mbox{ where }\theta:=\Arg(m).
\]
By elementary trigonometry, one can check that we then have
\[
\Gamma=\rme^{i\theta}+\frac{4}{\sqrt{3}}\sin\left(\frac{\abs{\theta}}
{3}\right)\rme^{2(\theta-\pi)i/3}
\]
or
\[
\Gamma=-\rme^{i\theta}+\frac{4}{\sqrt{3}}\sin\left(\frac{2\pi}{3}+
\frac{\abs{\theta}} {3}\right)\rme^{2\theta i/3}.
\]
Recalling our definition of the Airy functions $u_{11}(x)$ and
$u_{12}(x)$ (\ref{defu}), we put
\[
z(h,\lam,x):=h^{-2/3}\rme^{-2\theta i/3}(x\rme^{\theta i}-\lam),
\]
and can rewrite (\ref{asymptotics1}) explicitly in terms of $h$.
Then, taking the modulus we obtain
\[
\abs{Ai(z(h,\lam,x))}
\]
\begin{equation}\label{asyh}
=h^{1/6}\frac{\abs{x\rme^{\theta
i}-\lam}}{2\sqrt\pi}^{-1/4}\exp\left(-\frac{2}{3h}
\Re(\rme^{-2\theta i/3}(x\rme^{\theta i}-\lam))^{3/2}\right)
(1+O(h))
\end{equation}
as $h\to 0$, valid for $\abs{\Arg(z(h,\lam,x))}<\pi$. Therefore,
in order to estimate the moduli of the Airy functions
$Ai_k(z(h,\lam,x))$ in the limit $h\to 0$, it is sufficient to
examine the behaviour of the functions
\[ x\mapsto\Re\{(\rme^{-2k\pi i/3}z(h,\lam,x))^{3/2}\},\qquad
x\in\R
\]
for $k=-1,0,1.$

The basic idea of our proof is to show that for $\lam$ outside an
arbitrarily small $\eps$-neighbourhood of $Y(V(-1),V(1))$, one can
assign allowable values of $j$ and $k$ from $\{-1,0,1\}$ (in the
sense of (\ref{allowable})) such that the inequalities
\begin{equation}\label{eq1}
\Re\{(\rme^{-2j\pi i/3}z(h,\lam,-1)^{3/2}\}< \Re\{(\rme^{-2j\pi
i/3}z(h,\lam,1)^{3/2}\}
\end{equation}
and
\begin{equation}\label{eq2}
\Re\{(\rme^{-2k\pi i/3}z(h,\lam,1)^{3/2}\}< \Re\{(\rme^{-2k\pi
i/3}z(h,\lam,-1)^{3/2}\}
\end{equation}
hold in the limit $h\to 0$. This will then be enough, by
(\ref{asyh}), to ensure that $u_{11}(1)u_{12}(-1)$ and
$u_{11}(-1)u_{12}(1)$ are of different orders of magnitude as
$h\to 0$, thus implying (\ref{state1}).

Using the statements of the previous section; for all values of
$\lam$ such that
\[ z(h,\lam,\pm 1):=h^{-2/3}\rme^{-2\theta
i/3}(\pm\rme^{\theta i}-\lam) \] does not lie within an
$\eps$-neighbourhood of any of the Stokes' lines, and
$z(h,\lam,-1)$, $z(h,\lam,1)$ lie in different principal domains,
one can always obtain (\ref{eq1}) and (\ref{eq2}), and the
asymptotics (\ref{asyh}) will be valid. However, for $\lam$ lying
in the sector having its apex at $\Gamma$, and bounded by the rays
\[
\Gamma+r\rme^{2\theta i/3}\qquad\mbox{ and }\qquad
\left\{\begin{array}{ll} \Gamma+r\rme^{2\theta i/3+2\pi i/3}
&\mbox{ for }\theta<0 \\ \Gamma+r\rme^{2\theta i/3-2\pi i/3}
&\mbox{ for }\theta\geq 0 \\
\end{array} \right.
\]
$r\in[0,\infty)$, it is easy to check that $\rme^{-2k\pi
i/3}z(h,\lam,\pm 1)$ both lie in the same principal domain, for
each $k\in\{-1,0,1\}$. Then it is also straightforward to check
that as $x$ ranges from $-1$ to $1$, the function
\[
x\mapsto\Re\{(\rme^{-2k\pi i/3}z(\eps,\lam,x))^{3/2}
\]
which has been at the heart of our analysis, has a single
maximum/minimum. Together with the identities (\ref{id1}),
(\ref{id2}) and (\ref{id3}), this means that there $will$ be
values of $\lam$ such that equality holds in both (\ref{eq1}) and
(\ref{eq2}) - no matter what choices of $j$ and $k$ are made.
Thus, (and without loss of generality assuming $j=k=0$,) the set
of $\lam$ satisfying
\begin{equation}\label{gameq}
\Re\{(\rme^{-2\theta i/3}(\rme^{\theta i}-\lam))^{3/2}\}=
\Re\{(\rme^{-2\theta i/3}(-\rme^{\theta i}-\lam))^{3/2}\}
\end{equation}
lies in $Y(V(-1),V(1))$. We now examine this set in more detail.
Expanding the Taylor series, we have
\[(\rme^{-2\theta i/3}(\rme^{\theta
i}-\lam))^{3/2}=-i\lam^{3/2}\rme^{-\theta
i}+\frac{3}{2}i\lam^{1/2}-\frac{3}{8}i\lam^{-1/2}\rme^{\theta
i}-\frac{1}{16}i\lam^{-3/2}\rme^{2\theta i}+O(\lam^{-5/2})
\]
and
\[
(\rme^{-2\theta i/3}(-\rme^{\theta
i}-\lam))^{3/2}=-i\lam^{3/2}\rme^{-\theta
i}-\frac{3}{2}i\lam^{1/2}-\frac{3}{8}i\lam^{-1/2}\rme^{\theta
i}+\frac{1}{16}i\lam^{-3/2}\rme^{2\theta i}+O(\lam^{-5/2})
\]
as $\abs{\lam}\to\infty$. Dividing through by $-i$, this means
that (\ref{gameq}) will hold if and only if
\[
\Im\left\{\frac{3}{2}\lam^{1/2}-\frac{\lam^{-3/2}\rme^{2\theta
i}}{16}\right\}=O(\lam^{-5/2})
\]
as $\abs{\lam}\to\infty$. Putting $\lam:=\rho\rme^{\phi i}$, this
is equivalent to the requirement
\[
\sin\left(\frac{\phi}{2}\right)-\frac{1}{24\rho^2}
\sin\left(2\theta-\frac{3\phi}{2}\right)=O(\rho^{-6})
\]
as $\rho\to\infty$. But then
\begin{eqnarray*}
\Im(\lam)&&=\rho\sin(\phi)\\
         &&=2\rho\sin\left(\frac{\phi}{2}\right)\cos\left(\frac{\phi}{2}\right)\\
         &&=2\rho\cos\left(\frac{\phi}{2}\right)\left\{\frac{1}{24\rho^2}\sin
         \left(2\theta-\frac{3\theta}
         {2}\right)+O(\rho^{-6})\right\}\\
         &&=O(\rho^{-1})\\
\end{eqnarray*}
as $\rho\to\infty$. By our definition of $\Gamma$, $z(h,\Gamma,\pm
1)$ lies at the intersection of two Stokes' lines, and so
\[
\Re\{(z(h,\Gamma,\pm 1))^{3/2}\}=0
\]
showing that $\Gamma$ certainly lies in the set of $\lam$
satisfying (\ref{gameq}). Therefore, we deduce that
$Y(V(-1),V(1))$ contains a curve from $\Gamma$ asymptotic to the
positive real-axis.

Finally, we must examine what happens when $z$ does lie on one of
the Stokes' lines.

Firstly, suppose $\Arg(z(h,\lam,1))=\pi/3$, corresponding to
$\lam$ lying on the ray centred at $\rme^{\theta i}$ and passing
through $\Gamma$. Then $k=0,1$ are allowable, and one checks that
if $\lam$ lies on the segment $[\rme^{\theta i},\Gamma)$, we have
$z(h,\lam,-1)\in S_{-1}$. It follows by (\ref{id2}) that
\[
\Re\{(\rme^{-2\pi
i/3}z(h,\lam,-1))^{3/2}\}=\Re\{(z(h,\lam,-1))^{3/2}\},
\]
and so (\ref{eq1}) and (\ref{eq2}) cannot hold. However, if $\lam$
lies on that part of the ray which extends past $\Gamma$ (but not
$\lam=\Gamma$ itself,) then $z(h,\lam,-1)\in S_{1}$, and
\[
\Re\{(\rme^{-2\pi
i/3}z(h,\lam,-1))^{3/2}\}=-\Re\{(z(h,\lam,-1))^{3/2}\}
\]
causing (\ref{eq1}) and (\ref{eq2}) to hold for $j=0$, $k=1$.

An entirely similar argument holds when $\Arg(z(h,\lam,-1))=\pi$,
corresponding to $\lam$ lying on the ray centred at $-\rme^{\theta
i}$ and passing through $\Gamma$ using (\ref{id1}), with $j=1$ and
$k=-1$.

Finally, the case where $\Arg(z(h,\lam,\pm 1))=-\pi/3$ is taken
care of using (\ref{id3}), which shows that we may use allowable
values $-1$ and $0$ to obtain (\ref{eq1}) and (\ref{eq2}).

This completes the proof.

In the case $\theta=\pi/2$; $\Gamma=1/\sqrt{3}$ lies on the
real-axis, and the figure $Y(-i,i)$ has three $linear$ `arms'.
When $\theta=0$, the symmetric case, $Y(-1,1)$ is the
semi-infinite interval $[-1,\infty)$, as is well-known from the
theory of self-adjoint operators.

We now give our main result.

\begin{theorem}\label{t1}
Let
\[  H_hf(x):=-h^2\frac{\rmd^2f(x)}{\rmd x^2}+V(x)f(x) \]
act on $L^2(-1,1)$ with Dirichlet boundary conditions, where $h>0$
is small, and $V(x)$ is the complex valued n-times piecewise
linear function
\[ V(x):=\left\{\begin{array}{ll} m_1x+l_1 & x_0\leq x<x_1 \\
m_2x+l_2 & x_1<x<x_2 \\ \qquad \vdots & \qquad \vdots
\\ m_nx+l_n & x_{n-1}<x\leq x_n \\
\end{array} \right. \]
with $-1=x_0<x_1<\ldots<x_n=1$ and the $m_i$, $l_i$,
$i=1,\ldots,n$ complex constants. We assume for each $i$ that if
$m_ix_i+l_i=m_{i+1}x_i+l_{i+1}$, then $m_i\neq m_{i+1}$. Put
$\theta_i:=\Arg(m_i)$, and, using our earlier notation
\[
T:=\bigcup_{i=1}^nY(V(x_i),V(x_{i+1})).
\]
Let $\eps>0$ and $N\in\Z^+$ be given. Then
\[
\Spec(H_h)\cap\{z:\abs{z}\leq N\}\subset\Nhd(T;\eps)
\]
for all small enough $h>0$.
\end{theorem}
\Proof Our proof involves an analysis of the behaviour of the
characteristic-determinant, i.e. the left-hand side of
(\ref{matrix1}), as $h\to 0$. We give a proof for the case $n=3$;
the general case follows by a similar argument. For $n=3$, the
characteristic-determinant is given by
\begin{equation}\label{chardet1}
\abs{ \matrix{ u_{11}(-1) & u_{12}(-1) & 0 & 0 & 0 & 0 \cr
u_{11}(x_1) & u_{12}(x_1) & -u_{21}(x_1) & -u_{22}(x_1) & 0 & 0
\cr u'_{11}(x_1) & u'_{12}(x_1) & -u'_{21}(x_1) &-u'_{22}(x_1) & 0
& 0 \cr 0 & 0 & u_{21}(x_2) & u_{22}(x_2) & -u_{31}(x_2) &
-u_{32}(x_2) \cr 0 & 0 & u'_{21}(x_2) & u'_{22}(x_2) &
-u'_{31}(x_2) & -u'_{32}(x_2) \cr 0 & 0 & 0 & 0 & u_{31}(1) &
u_{32}(1) \cr }}
\end{equation}
and we must prove that for certain values of $\lam\in\C$, this
determinant does not vanish as $h\to 0$. Expanding
(\ref{chardet1}), one obtains
\[
\Big\{(u_{11}(-1)u_{12}(x_1)-u_{12}(-1)u_{11}(x_1))(u'_{22}(x_1)u_{21}(x_2)-
u'_{21}(x_1)u_{22}(x_2))\times \]
\[ \times(u_{31}(1)u'_{32}(x_2)-u'_{31}(x_2)u_{32}(1))\Big\}-
\]
\[
-\Big\{(u_{11}(-1)u_{12}(x_1)-u_{12}(-1)u_{11}(x_1))(u'_{22}(x_1)u'_{21}(x_2)-
u'_{21}(x_1)u'_{22}(x_2))\times \]
\[ \times(u_{31}(1)u_{32}(x_2)-u_{31}(x_2)u_{32}(1))\Big\}+
\]
\[
+\Big\{(u_{11}(-1)u'_{12}(x_1)-u_{12}(-1)u'_{11}(x_1))(u_{22}(x_1)u'_{21}(x_2)-
u_{21}(x_1)u'_{22}(x_2))\times \]
\[ \times(u_{31}(1)u_{32}(x_2)-u_{31}(x_2)u_{32}(1))\Big\}-
\]
\[ -\Big\{(u_{11}(-1)u'_{12}(x_1)-u_{12}(-1)u'_{11}(x_1))(u_{22}(x_1)u_{21}(x_2)-
u_{21}(x_1)u_{22}(x_2))\times \]
\begin{equation}\label{chardet2}
\times(u_{31}(1)u'_{32}(x_2)-u'_{31}(x_2)u_{32}(1))\Big\}
\end{equation}
where so far, no asymptotics are involved.

Now, let $\eps>0$ and $N\in\Z^+$ be as given in the statement of
the theorem. Taking any
\[ \lam\in\{z:\abs{z}\leq N\}\backslash\Nhd(T;\eps), \]
we can use the results of Proposition \ref{t2} to show that
(\ref{chardet1}) is non-zero in the limit as $h\to 0$. Indeed, by
the proof of Proposition \ref{t2}, we can ensure that the
asymptotic estimates
\[
u_{12}(-1)u_{11}(x_1)=o(u_{11}(-1)u_{12}(x_1)), \]
\[ u_{21}(x_1)u_{22}(x_2)=o(u_{22}(x_1)u_{21}(x_2)) \]and
\[ u_{31}(x_2)u_{32}(1)=o(u_{32}(x_2)u_{31}(1))
\]
hold, as $h\to 0$. Then, using the standard asymptotic expansions
of the Airy functions \cite{Olver74}, which give
\begin{equation}\label{asymptotics}
Ai(z)=\frac{z^{-1/4}}{2\sqrt\pi}\exp\left(-\frac{2}{3}z^{3/2}\right)
(1+O(z^{-3/2}))
\end{equation}
and
\begin{equation}\label{asymptotics'}
Ai'(z)=-\frac{z^{1/4}}{2\sqrt\pi}\exp\left(-\frac{2}{3}z^{3/2}\right)
(1+O(z^{-3/2}))
\end{equation}
as $\abs{z}\to\infty$, valid for all $z$ such that
$\abs{\Arg(z)}<\pi$; we see that, if
\[
z(h,\lam,x_i):=h^{-2/3}m_i^{-2/3}((m_ix_i+l_i)-\lam),
\]
then
\begin{eqnarray*}
\frac{\rmd}{\rmd x}Ai(z)&&=\frac{\rmd z}{\rmd x}Ai'(z)
\\ &&=-\frac{\rmd z}{\rmd x}\frac{
z^{1/4}}{2\sqrt\pi}\exp\left(-\frac{2}{3}z^{3/2}\right)(1+O(z^{-3/2}))\\
&&=-\frac{\rmd z}{\rmd x}z^{1/2}
\frac{z^{-1/4}}{2\sqrt\pi}\exp\left(-\frac{2}{3}z^{3/2}\right)(1+O(z^{-3/2}))\\
\end{eqnarray*}
as $\abs{z}\to\infty$. Comparing this last expression with
(\ref{asymptotics}), and using $f\sim g$ to mean that
\[
\frac{f(h)}{g(h)}\to 1 \qquad\mbox{ as }h\to 0,
\]
we obtain
\[
\frac{\rmd}{\rmd
x}Ai(z(h,\lam,x_i))\sim-h^{-1}((m_ix_i+l_i)-\lam)^{1/2}Ai(z(h,\lam,x_i))
\qquad\mbox{ as }h\to 0.
\]
Moreover, similar calculations show that
\[
\frac{\rmd}{\rmd x}Ai(\rme^{\pm 2\pi i/3}z(h,\lam,x_i)) \sim
h^{-1}((m_ix_i+l_i)-\lam))^{1/2}Ai(z(h,\lam,x_i))\qquad\mbox{ as
}h\to 0.
\]
Reverting to our notation of (\ref{defu}), we will write
\begin{equation}\label{sim1}
u'_{i1}(x_i)\sim h^{-1}c_{i1}(x_i)u_{i1}(x_i)\qquad \mbox{ etc.}
\end{equation}
as $h\to 0$, where it is important to note that the $c_{ij}(x_i)$,
$i=1,\ldots,(n-1)$, $j=1,2$ are independent of $h$. Then, since
the constant terms $c_{ij}(x_i)$ are negligible in magnitude
compared with the exponential terms $u_{ij}(x_i)$ as $h\to 0$, the
relations (\ref{sim1}) imply that we also have the estimates
\[
u_{12}(-1)u'_{11}(x_1)=o(u_{11}(-1)u'_{12}(x_1)),
\]
\[ u'_{21}(x_1)u_{22}(x_2)=o(u'_{22}(x_1)u_{21}(x_2)) \]
\[ u'_{31}(x_2)u_{32}(1)=o(u'_{32}(x_2)u_{31}(1))
\]
\[ u'_{21}(x_1)u'_{22}(x_2)=o(u'_{22}(x_1)u'_{21}(x_2))
\]and
\[ u_{21}(x_1)u'_{22}(x_2)=o(u_{22}(x_1)u'_{21}(x_2))
\]
as $h\to 0$. Returning to (\ref{chardet2}), we first use the above
estimates (since we may ignore the sub-dominant term in each
round-bracketed expression), and then the relations (\ref{sim1})
again, to obtain the asymptotic estimate on the first of the
curly-bracketed terms:
\[
\Big\{(u_{11}(-1)u_{12}(x_1)-u_{12}(-1)u_{11}(x_1))(u'_{22}(x_1)u_{21}(x_2)-
u'_{21}(x_1)u_{22}(x_2))\times \]
\[ \times(u_{31}(1)u'_{32}(x_2)-u'_{31}(x_2)u_{32}(1))\Big\}
\]
\begin{eqnarray*}
&&\sim
u_{11}(-1)u_{12}(x_1)u'_{22}(x_1)u_{21}(x_2)u_{31}(1)u'_{32}(x_2)\\
&&\sim u_{11}(-1)u_{12}(x_1)\eps^{-1/2}c_{22}(x_1)u_{22}(x_1)
u_{21}(x_2)u_{31}(1)\eps^{-1/2}c_{32}(x_2)u_{32}(x_2)\\
&&=h^{-2}[c_{22}(x_1)c_{32}(x_2)](u_{11}(-1)u_{12}(x_1)u_{22}(x_1)
u_{21}(x_2)u_{31}(1)u_{32}(x_2))\\
\end{eqnarray*}
as $h\to 0$. Similar estimates apply to each of the remaining
three terms in (\ref{chardet2}) i.e.
\[
\Big\{(u_{11}(-1)u_{12}(x_1)-u_{12}(-1)u_{11}(x_1))(u'_{22}(x_1)u'_{21}(x_2)-
u'_{21}(x_1)u'_{22}(x_2))\times \]
\[ \times(u_{31}(1)u_{32}(x_2)-u_{31}(x_2)u_{32}(1))\Big\}
\]
\[
\sim
h^{-2}[c_{22}(x_1)c_{21}(x_2)](u_{11}(-1)u_{12}(x_1)u_{22}(x_1)
u_{21}(x_2)u_{31}(1)u_{32}(x_2)),
\]
\par
\[
\Big\{(u_{11}(-1)u'_{12}(x_1)-u_{12}(-1)u'_{11}(x_1))(u_{22}(x_1)u'_{21}(x_2)-
u_{21}(x_1)u'_{22}(x_2))\times \]
\[ \times(u_{31}(1)u_{32}(x_2)-u_{31}(x_2)u_{32}(1))\Big\}
\]
\[
\sim
h^{-2}[c_{12}(x_1)c_{21}(x_2)](u_{11}(-1)u_{12}(x_1)u_{22}(x_1)
u_{21}(x_2)u_{31}(1)u_{32}(x_2))
\]and
\[
\Big\{(u_{11}(-1)u'_{12}(x_1)-u_{12}(-1)u'_{11}(x_1))(u_{22}(x_1)u_{21}(x_2)-
u_{21}(x_1)u_{22}(x_2))\times \]
\[ \times(u_{31}(1)u'_{32}(x_2)-u'_{31}(x_2)u_{32}(1))\Big\}
\]
\[
\sim
h^{-2}[c_{12}(x_1)c_{32}(x_2)](u_{11}(-1)u_{12}(x_1)u_{22}(x_1)
u_{21}(x_2)u_{31}(1)u_{32}(x_2))
\]
as $h\to 0$. Collecting these estimates together, we see that the
characteristic determinant (\ref{chardet1}) tends asymptotically
towards
\[
h^{-2}\{(c_{22}(x_1)-c_{12}(x_1))(c_{32}(x_2)-c_{21}(x_2))\}
(u_{11}(-1)u_{12}(x_1)u_{22}(x_1) u_{21}(x_2)u_{31}(1)u_{32}(x_2))
\]
as $h\to 0$. The Airy functions $Ai(z)$ have countably many
negative real zeros, \cite{Olver74}; and so by our choice of
$\lam$ outside $\Nhd(T;\eps)$ together with the proof of
Proposition \ref{t2}, we are assured that none of the Airy
functions $u_{ij}(x_i)$ vanishes. Therefore, the determinant
(\ref{chardet1}) does not vanish in the limit as $h\to 0$,
provided the `constant' terms
\begin{equation}\label{eqncon}
c_{22}(x_1)\neq c_{12}(x_1)\qquad\mbox{ and }\qquad
c_{32}(x_2)\neq c_{21}(x_2).
\end{equation}
Our choice of $\lam$ ensures that each of the individual constant
terms $c_{ij}(x_i)$ is non-zero. Moreover, reviewing the proof of
Proposition \ref{t2} and the identities (\ref{id1})-(\ref{id3}),
we see that the choices for $j$ and $k$ are not uniquely
determined. Therefore, it is always possible to ensure that
(\ref{eqncon}) holds, even when $V$ is continuous at some or all
of the $x_is$. For example, if it happens that $V(x_1-)=V(x_1+)$,
then we choose $j$ and $k$ so that the constants $c_{12}(x_1)$ and
$c_{22}(x_1)$ take different signs (by the calculations
immediately above (\ref{sim1})). Thus, we deduce that such $\lam$
cannot be an eigenvalue.

It now just requires the following compactness argument to
complete the proof. Let $B(z;\eps)$ denote the open ball centred
at $z$, with radius $\eps$. Our argument so far shows that for any
\[ \lam\in\{z\in\C:\abs{z}\leq N\} \] such that
\[ B(\lam;\eps)\cap T=\emptyset \]
we have
\[ B(\lam;\eps)\cap\Spec(H_h)=\emptyset \]
for all $0<h<E_{\lam}$, where $E_{\lam}$ is some positive constant
dependent upon $\lam$. Let
\[
M:=\{z\in\C:\abs{z}\leq N,\mbox{ and }\dist(z,T)\geq 2\eps\},
\]
so that $M$ is compact. Then for all $\lam\in M$
\[
B(\lam;\eps)\cap\Nhd(T;\eps)=\emptyset
\]
and so
\[
M\subseteq \bigcup_{\lam\in M}B(\lam;\eps).
\]
But by compactness this means that there exists a finite
sub-covering
\[ M\subseteq \bigcup_{r=1}^n B(\lam_r;\eps_{\lam_r}).  \]
Taking $E$ to be $\min(E_{\lam_1},\dots ,E_{\lam_n})>0$, we deduce
that for all $0<h<E$ we have
\[ \Spec(H_h)\cap M=\emptyset \]
and this is equivalent to the statement of the theorem.

\begin{remark}
An important but subtle point, to note is that the zeros of
\begin{equation}\label{zeros}
(u_{11}(-1)u_{12}(x_1)u_{22}(x_1) u_{21}(x_2)u_{31}(1)u_{32}(x_2))
\end{equation}
as a function of $\lam$, are $not$ the same as the zeros of
(\ref{chardet1}). However, by a similar argument to that of
Shkalikov \cite{Shk97} (i.e. using (\ref{asymptotics})), one can
readily show that along each of the bounded arms of the Y-shaped
figures making up $T$, the zeros (eigenvalues) do converge as
$h\to 0$ to form a dense set. Finding an asymptotic expression for
the density of spectral points along the infinite lines (in the
direction of the positive real-axis) appears to be a much more
difficult problem; and we have no results yet in that direction.
\end{remark}

To illustrate Theorem \ref{t1}, in figure 3 we show a Matlab plot
of the discretised version of the problem where the potential
\[ V(x):=\left\{\begin{array}{ll}
2ix+i &\mbox{ for $-1\leq x<0$ }\\ (i+1)x &\mbox{ for $0<x\leq 1
$}\\
\end{array} \right. \]

We now return to the operator $H_{\del,h}$ defined in the first
section.
\begin{corollary}\label{c1}
Let $H_{\del,h}$ be the non-self-adjoint operator defined by
\[ H_{\del,h}:=-h^2\frac{\rmd^2}{\rmd x^2}+V_{\del}(x) \]
acting on $L^2(-1,1)$ with Dirichlet boundary conditions, $h>0$ ,
and
\[ V_{\del}(x):=\left\{\begin{array}{ll}
i(x-\del) &\mbox{ for $x<0$ }\\ i(x+\del) &\mbox{ for $x>0$}\\
\end{array} \right. \] with $\del>0$. Define
$S\subset\C$ to be the double Y-shaped figure given by the line
segments
\[ [i\del,1/2\sqrt{3}+i(1+2\del)/2] \]
\[ [i(1+\del),1/2\sqrt{3}+i(1+2\del)/2] \]together with
\[ [1/2\sqrt{3}+i(1+2\del)/2,+\infty), \]
and
\[ [-i\del,1/2\sqrt{3}-i(1+2\del)/2] \]
\[ [-i(1+\del),1/2\sqrt{3}-i(1+2\del)/2] \]together with
\[ [1/2\sqrt{3}-i(1+2\del)/2,+\infty). \]
Then, given any $\eps>0$ and $N\in\Z^+$, we have
\[
\Spec(H_{\del,h})\cap\{z:\abs{z}\leq n\}\subset\Nhd(S;\eps)
\]
for small enough $h>0$ (see figure 1).
\end{corollary}
By analyticity, however, for fixed $h>0$ we have
\[ \lim_{\del\to 0}\Spec(H_{\del,h})=\Spec(H_{0,h}). \]
Hence, $\lim_{h\to 0}\lim_{\del\to 0}\Spec(H_{\del,h})$ is
contained within an arbitrarily small neighbourhood of the line
segments
\[ [i,1/\sqrt{3}],[-i,1/\sqrt{3}]\qquad\mbox{ and }\qquad
[1/\sqrt{3},\infty) \] (see figure 2). Thus, the operations of
taking limits do not commute, in the sense that as sets
\[ \lim_{h\to 0}\;\lim_{\del\to 0}\;\Spec(H_{\del,h})\neq
\lim_{\del\to 0}\;\lim_{h\to 0}\;\Spec(H_{\del,h}). \]


\section{Simultaneous limits for $H_{\del,h}$}\label{s3}

In response to questions which we have been asked, we give an
analysis for the situation in which $\del$ and $h$ of Corollary
\ref{c1} are not independent of each other.
\begin{proposition}\label{t5}
Defining the operator $H_{\del,h}$ as above, and putting
\[ \del:=h^{1/p} \]
we have
\[
\lim_{h\to 0}\Spec(H_{h^{1/p},h})=\lim_{h\to
0}\Spec(H_{0,h})\qquad\mbox{ if }0<p<1
\]
and
\[
\lim_{h\to 0}\Spec(H_{h^{1/p},h})= \lim_{\del\to 0}\;\lim_{h\to
0}\;\Spec(H_{\del,h})\qquad\mbox{ if }p\geq 1.
\]
\end{proposition}
\Proof Referring to (\ref{matrix1}) and expanding, we see that the
characteristic determinant of $H_{\del,h}$ is given by
\[ \Big\{(u_{11}(-1)u_{12}(0)-u_{12}(-1)u_{11}(0))
(u'_{22}(0)u_{21}(1)-u'_{21}(0)u_{22}(1))\Big\}- \]
\begin{equation}\label{chardet3}
-\Big\{(u_{11}(-1)u'_{12}(0)-u_{12}(-1)u'_{11}(0))
(u_{22}(0)u_{21}(1)-u_{21}(0)u_{22}(1))\Big\},
\end{equation}
whereas the characteristic determinant of $H_{0,h}$ is given by
\begin{equation}\label{chardet4}
u_{21}(1)u_{12}(-1)-u_{22}(1)u_{11}(-1).
\end{equation}
Now, putting

$u_{12}(0):=Ai_k(h^{-2/3}\rme^{-\pi i/3}(-i\del-\lam))$  and
$u_{22}(0):=Ai_k(h^{-2/3}\rme^{-\pi i/3}(i\del-\lam))$

it is clear by analyticity, that
\begin{equation}\label{est1}
u_{12}(0)\sim u_{22}(0)\qquad\mbox{ and }\qquad u_{11}(0)\sim
u_{21}(0)
\end{equation}
as $\del\to 0$. Moreover, the calculations preceding (\ref{sim1})
show that
\begin{equation}\label{est2}
u_{i1}'(0)\sim -\del^{-p}(-\lam)^{1/2}u_{i1}(0)\qquad\mbox{ and
}\qquad u_{i2}'(0)\sim \del^{-p}(-\lam)^{1/2}u_{i2}(0)
\end{equation}
as $\del\to 0$, for $i=1,2$. Therefore, using first (\ref{est2})
and then (\ref{est1}), the characteristic determinant
(\ref{chardet3}) tends asymptotically towards
\[
2\del^{-p}(-\lam)^{1/2}\Big(u_{11}(-1)u_{12}(0)u_{21}(0)u_{22}(1)-
u_{12}(-1)u_{11}(0)u_{22}(0)u_{21}(1)\Big)
\]
\[
\sim
2\del^{-p}(-\lam)^{1/2}\Big(u_{22}(1)u_{11}(-1)-u_{21}(1)u_{12}(-1)\Big)
\]
as $\del\to 0$. Then the zeros of (\ref{chardet3}) tend
asymptotically toward the zeros of (\ref{chardet4}) by $\rm
Rouch\acute e$'s theorem, explaining the behaviour of
\[
\lim_{h\to 0}\;\lim_{\del\to 0}\;\Spec(H_{\del,h}).
\]

Substituting $\del=h^{1/p}$, the character of $\lim_{h\to
0}\Spec(H_{h^{1/p},h})$ therefore depends upon the range of $p$
for which
\[
\frac{u_{22}(0)}{u_{12}(0)}\to 1\qquad\mbox{ and }\qquad
\frac{u_{21}(0)}{u_{11}(0)}\to 1
\]
as $h\to 0$. Now, without loss of generality, and using our
earlier notation, let
\[ \frac{u_{22}(0)}{u_{12}(0)}:=\frac{Ai_{-1}(z_1)}{Ai_{-1}(z_2)} \]
where
\[ z_1:=h^{-2/3}\rme ^{\pi i/3}(ih^{1/p}-\lam)\qquad\mbox{
and }\qquad z_2:=h^{-2/3}\rme ^{\pi i/3}(-ih^{1/p}-\lam) \] so
that, using the standard asymptotics (\ref{asymptotics1})
\begin{eqnarray*}
\frac{u_{22}(0)}{u_{12}(0)}&&=\frac{z_1^{-1/4}\exp\left
(-\frac{2}{3}z_1^{3/2}\right)(1+O(z_1^{-3/2}))}{z_2^{-1/4}\exp\left
(-\frac{2}{3}z_2^{3/2}\right)(1+O(z_2^{-3/2}))}\\
&&=\left(\frac{z_1}{z_2}\right)^{-1/4}\exp\left(-\frac{2}{3}
\left[z_1^{3/2}-z_2^{3/2}\right]\right) (1+O(z_1^{-3/2}))\\
&&\sim\exp\left(-\frac{2}{3}\left[z_1^{3/2}-z_2^{3/2}\right]\right)\\
\end{eqnarray*}
as $h\to 0$. But
\begin{eqnarray*}
z_1^{3/2}-z_2^{3/2}&&=h^{-1}\rme^{\pi
i/3}\left\{(ih^{1/p}-\lam)^{3/2}-(-ih^{1/p}-\lam)^{3/2}\right\}\\
&&=h^{-1}\rme^{\pi
i/3}(-\lam)^{3/2}\left\{(1-ih^{1/p}/\lam)^{3/2}-(1+ih^{1/p}/\lam)^{3/2}\right\}\\
&&=h^{-1}\rme^{\pi
i/3}(-\lam)^{3/2}\left\{(1-3ih^{1/p}/2\lam+\dots)
-(1+3ih^{1/p}/2\lam+\dots)\right\}\\ &&=h^{-1}\rme^{\pi
i/3}(-\lam)^{3/2}(-3ih^{1/p}/\lam+\dots)\\ &&\to 0\\
\end{eqnarray*}
as $h\to 0$ if and only if $0<p<1$. So, provided $0<p<1$
\[ \frac{u_{22}(0)}{u_{12}(0)}\to 1 \qquad\mbox{ as }h\to 0 \]
and a similar calculation shows that we then also have
\[ \frac{u_{21}(0)}{u_{11}(0)}\to 1 \qquad\mbox{ as }h\to 0, \]
as required.


\vskip 1in
{\bf Acknowledgements }I would like to thank E. B. Davies for
suggesting this problem, and his guidance in the course of solving
it.
\par

\vskip 0.3in Department of Mathematics \newline King's College
\newline Strand \newline London WC2R 2LS \newline England \\
e-mail:Redparth@mth.kcl.ac.uk \vfil
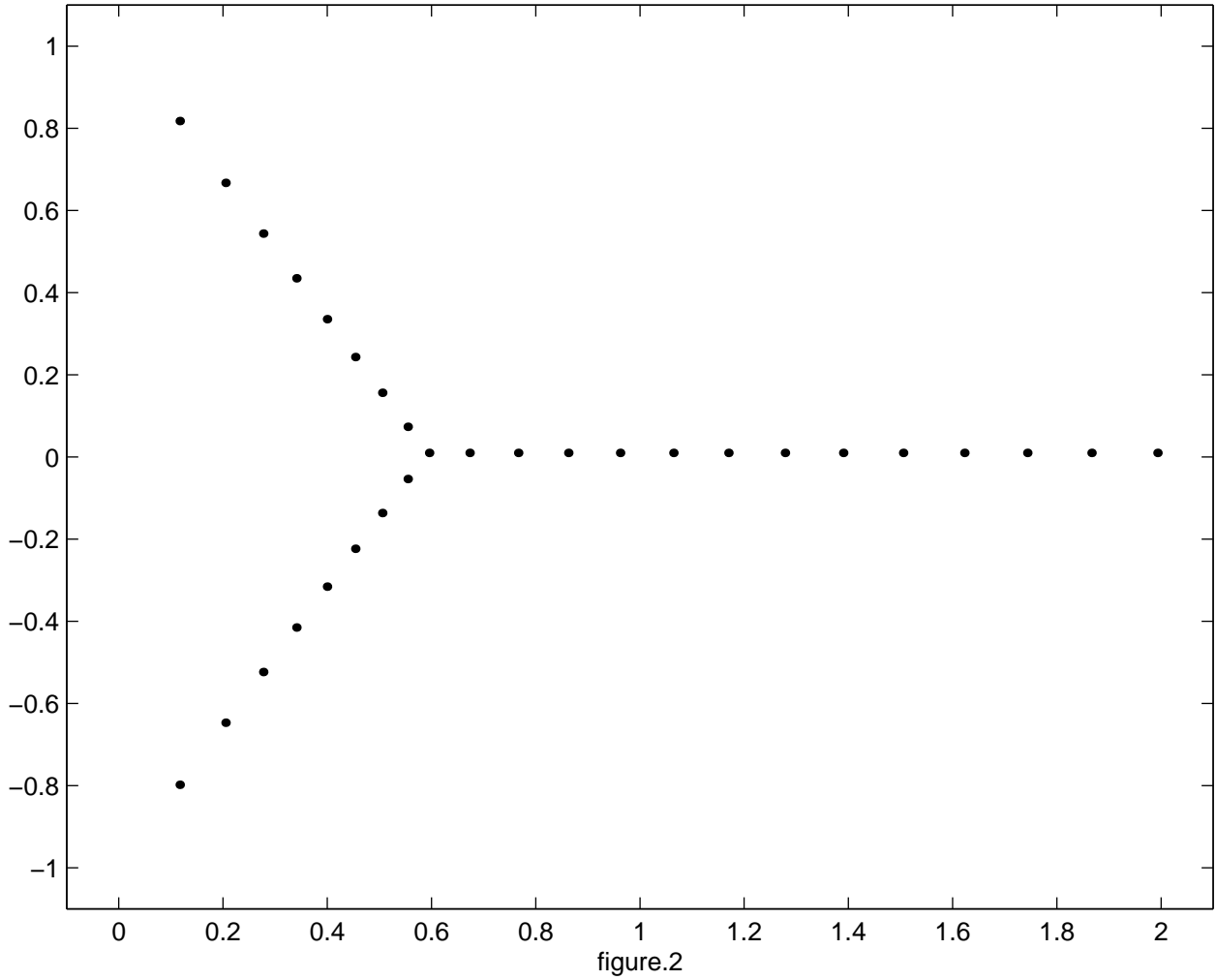
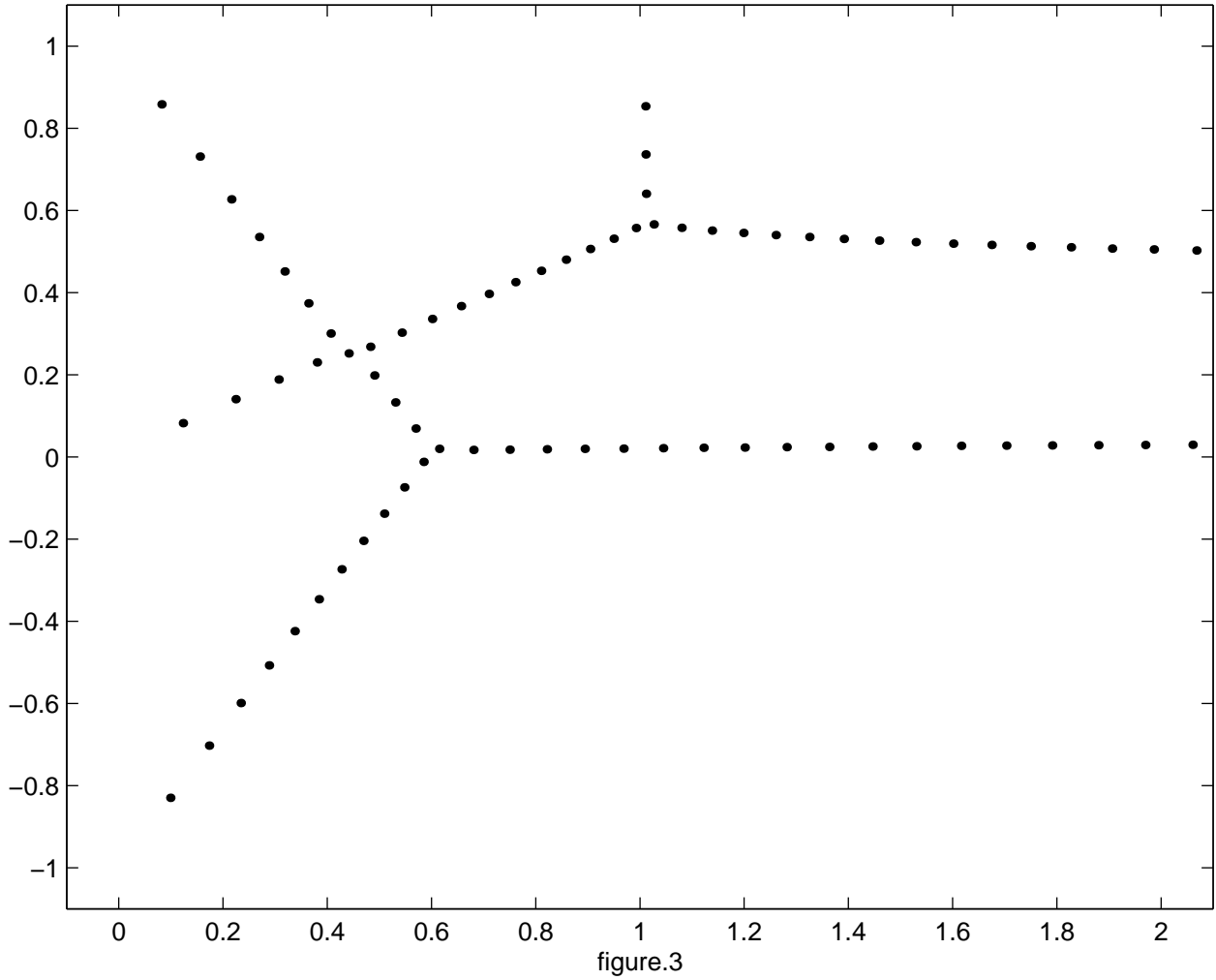
\end{document}